\documentclass[11pt,a4paper]{article}
\usepackage{amssymb}
\usepackage{amsmath}
\usepackage{amsthm}
\usepackage{verbatim}
\usepackage{pgf}
\usepackage{listings}
\usepackage{graphicx}
\usepackage{tabularx}
\usepackage{psfrag}
\usepackage{epstopdf}
\usepackage{multirow}
\usepackage{ifthen}
\usepackage{booktabs}
\usepackage{mathtools}
\usepackage{algorithm}
\usepackage{algorithmicx}
\usepackage{algpseudocode}
\usepackage{authblk}
\usepackage[justification=raggedleft,singlelinecheck=false]{caption}
\usepackage{geometry}
\allowdisplaybreaks

\numberwithin{equation}{section}

\newtheorem{thm}{Theorem}[section]  

\theoremstyle{definition}

\newtheorem{lmm}[thm]{Lemma}
\newtheorem{prp}[thm]{Proposition}
\newtheorem{trm}[thm]{Theorem}
\newtheorem{crl}[thm]{Corollary}

\newtheorem{rem}[thm]{Remark}

\theoremstyle{remark}
\newtheorem*{prf}{Proof}

\begin{document}

\title{Recurrence Relations for Values of the Riemann Zeta Function in Odd Integers}

\author{Tobias Kyrion}

\affil{H\"urth, Germany}

\renewcommand\Authands{ and }

\maketitle

\begin{abstract}
It is commonly known that $\zeta(2k) = q_{k}\frac{\zeta(2k + 2)}{\pi^2}$ with known rational numbers $q_{k}$. In this work we construct recurrence relations of the form $\sum_{k = 1}^{\infty}r_{k}\frac{\zeta(2k + 1)}{\pi^{2k}} = 0$ and show that series representations for the coefficients $r_{k} \in \mathbb{R}$ can be computed explicitly.
\end{abstract}

\tableofcontents

\section{Summary}

In the first section we show that $\frac{\cosh(x)}{\sinh(x)^{2N + 1}}$ can be expressed as linear combination of $\frac{\sinh(x)}{\cosh(x)^{2N + 1}}$ and $\frac{\cosh(2x)}{\sinh(2x)^{2k + 1}}$ for some $k \leq N$. We achieve this by proving four identities between certain rational functions. Then we show that the $2n$-th derivative of $\coth(x)$ can be expressed as linear combination of $\frac{\cosh(x)}{\sinh(x)^{2k + 1}}$ with $k$ ranging from $1$ to $n$. We prove some useful recurrence relations between the coefficients of the $\frac{\cosh(x)}{\sinh(x)^{2k + 1}}$'s and compute explicitly the inverse of the matrix formed by these coefficients. We derive our main result Theorem 3.3 - the limit identity for $\lim_{\alpha \to 0_{+}}\sum_{n = 1}^{\infty}\frac{1}{n}\frac{\sinh(\alpha n)}{\cosh(\alpha n)^{2N + 1}}$ for a fixed $N \in \mathbb{N}$ - by applying our previous findings on Ramanujan's famous identity for the Riemann zeta function values in odd integers. As an application we finally determine recurrence relations of the form $\sum_{k = 1}^{\infty}r_{k}\frac{\zeta(2k + 1)}{\pi^{2k}} = 0$.

\section{Preliminaries}

\rem{Throughout this work we set $\binom{n}{k} =  0$ for $k < 0$ and for $n < k$.}

\subsection{Identities for $\frac{\cosh(x)}{\sinh(x)^{2N + 1}}$}

\begin{prp}
Let $M \in \mathbb{N}$ and $z \in \mathbb{C}$. Then we have the following four relations between rational functions:
\begin{align*}
\frac{1}{2}\left(z + \frac{1}{z}\right)^{4M + 2} - \frac{1}{2}\left(z - \frac{1}{z}\right)^{4M + 2} & = \; \sum_{k = 0}^{M}2^{4k}\frac{4M + 2}{2k + 1}\binom{M + k}{2k}\left(z^{2} - \frac{1}{z^{2}}\right)^{2(M - k)} \\
\frac{1}{2}\left(z + \frac{1}{z}\right)^{4M} - \frac{1}{2}\left(z - \frac{1}{z}\right)^{4M} & = \; 4\left(z^{2} + \frac{1}{z^{2}}\right)\sum_{k = 0}^{M - 1}2^{4k}\frac{M - k}{2k + 1}\binom{M + k}{2k}\left(z^{2} - \frac{1}{z^{2}}\right)^{2(M - 1 - k)} \\
\frac{1}{2}\left(z + \frac{1}{z}\right)^{4M} + \frac{1}{2}\left(z - \frac{1}{z}\right)^{4M} & = \; \sum_{k = 0}^{M}2^{4k}\frac{M}{M + k}\binom{M + k}{2k}\left(z^{2} - \frac{1}{z^{2}}\right)^{2(M - k)} \\
\frac{1}{2}\left(z + \frac{1}{z}\right)^{4M + 2} + \frac{1}{2}\left(z - \frac{1}{z}\right)^{4M + 2} & = \; \left(z^{2} + \frac{1}{z^{2}}\right)\sum_{k = 0}^{M}2^{4k}\binom{M + k}{2k}\left(z^{2} - \frac{1}{z^{2}}\right)^{2(M - k)}.
\end{align*}
\end{prp}

\begin{prf}
After factoring out $z^{4(M - j)} + z^{4(j - M)}$ in both sides of the first relation and $z^{4(M - j) - 2} + z^{4(j - M) + 2}$ in the second, $z^{4(M - j)} + z^{4(j - M)}$ in the third and $z^{4(M - j) + 2} + z^{4(j - M) - 2}$ in the fourth we are by comparing coefficients left to show the following relations: 
\begin{align*}
\sum_{k = 0}^{j}(-1)^{j - k}2^{4k}\left(\binom{2(M - k)}{j - k} - \binom{2(M - k)}{j - 1 - k}\right)\binom{M + k}{2k} & = \binom{4M + 2}{2j}, &  j & = 0, ..., M \\
\sum_{k = 0}^{j}(-1)^{j - k}2^{4k}\frac{M}{M + k}\binom{2(M - k)}{j - k}\binom{M + k}{2k} & = \binom{4M}{2j}, & j & = 0, ..., M \\
\sum_{k = 0}^{j}(-1)^{j - k}2^{4k}\frac{4M + 2}{2k + 1}\binom{2(M - k)}{j - k}\binom{M + k}{2k} & = \binom{4M + 2}{2j + 1}, & j & = 0, ..., M \\
\sum_{k = 0}^{j}(-1)^{j - k}2^{4k + 2}\frac{M - k}{2k + 1}\left(\binom{2(M - 1 - k)}{j - k} - \binom{2(M - 1 - k)}{j - 1 - k}\right)\binom{M + k}{2k} & = \binom{4M}{2j + 1}, & j & = 0, ..., M - 1. 
\end{align*}

We denote with $S_{2j}^{4M + 2}$, $S_{2j}^{4M}$, $S_{2j + 1}^{4M + 2}$ and $S_{2j + 1}^{4M}$ the left hand sides of above relations. With some algebra we can establish 
\begin{align*}
(2M - 2j + 3)(4M + 2)jS_{2j}^{4M + 2} + 4(2M - 2j + 2)&(2M - 2j + 1)(3M - j + 2)S_{2(j - 1) + 1}^{4M + 2} \\
- \;\;2(2M - 2j + 2)(4M + 1)(4M + 2)S_{2(j - 1) + 1}^{4M} & = -(2M - 2j + 1)(2M - j + 2)(4M + 2)S_{2(j - 1)}^{4M + 2} \\
\\
jS_{2j}^{4M} - 2MS_{2(j - 1) + 1}^{4M} & = - (2M - j + 1)S_{2(j - 1)}^{4M} \\
\\
(2j + 1)S_{2j + 1}^{4M + 2} - (4M + 2)S_{2j}^{4M + 2} & = - (4M - 2j + 3)S_{2(j - 1) + 1}^{4M + 2} \\
\\
(2M + 1)S_{2j + 1}^{4M} - (2M - 2j)S_{2j + 1}^{4M + 2} & = (2M + 1)S_{2(j - 1) + 1}^{4M}.
\end{align*}

A computation shows that these relations are also fulfilled by the corresponding right hand sides of the top most equations. Thus the proof follows by induction over $j$.

\hfill $\square$
\end{prf}

Setting $z = \exp(x)$ for $x \in \mathbb{C}$ and some algebraic manipulations give
\crl{For $M \in \mathbb{N}$ and $x \in \mathbb{C}$ we have
\begin{equation}
\begin{split}
\frac{\cosh(x)}{\sinh(x)^{4M - 1}} = \;\; & - \frac{\sinh(x)}{\cosh(x)^{4M - 1}} + \sum_{k = 0}^{M}2^{2(M + k)}\frac{M}{M + k}\binom{M + k}{2k}\frac{1}{\sinh(2x)^{2(M + k) - 1}} \\
= \;\; & \frac{\sinh(x)}{\cosh(x)^{4M - 1}} + \sum_{k = 0}^{M - 1}2^{2(M + k) + 1}\frac{M - k}{2k + 1}\binom{M + k}{2k}\frac{\cosh(2x)}{\sinh(2x)^{2(M + k) + 1}} \label{eq_coth_odd}
\end{split}
\end{equation}
and
\begin{equation}
\begin{split}
\frac{\cosh(x)}{\sinh(x)^{4M + 1}} = \;\; & \frac{\sinh(x)}{\cosh(x)^{4M + 1}} + \sum_{k = 0}^{M}2^{2(M + k) + 1}\frac{2M + 1}{2k + 1}\binom{M + k}{2k}\frac{1}{\sinh(2x)^{2(M + k) + 1}} \\
= \;\; & - \frac{\sinh(x)}{\cosh(x)^{4M + 1}} + \sum_{k = 0}^{M}2^{2(M + k) + 1}\binom{M + k}{2k}\frac{\cosh(2x)}{\sinh(2x)^{2(M + k) + 1}}. \label{eq_coth_even}
\end{split}
\end{equation}
}

\subsection{The $2n$-th Derivative of $\coth(x)$ and Some Useful Recurrence and Matrix Relations}

The $2n$-th derivative of $\coth(x)$ will be of special interest in the next section.

\lmm{let $n \in \mathbb{N}$ and $x \in \mathbb{C}$. Then we have
\begin{align*}
\frac{d^{2n}}{dx^{2n}}\coth(x) = \sum_{k = 1}^{n}\frac{2}{4^{k}}\sum_{j = 1}^{k}(-1)^{k - j}\binom{2k}{k - j}(2j)^{2n}\frac{\cosh(x)}{\sinh(x)^{2k + 1}}.
\end{align*}
}

\begin{prf}
We make an ansatz of the form
\begin{align*}
\frac{d^{2n}}{dx^{2n}}\coth(x) = \sum_{k = 1}^{n}c_{n, k}\frac{\cosh(x)}{\sinh(x)^{2k + 1}},
\end{align*}
with coefficients $c_{n, k} \in \mathbb{R}$. Utilizing
\begin{align*}
\frac{d^{2}}{dx^{2}}\frac{\cosh(x)}{\sinh(x)^{2k + 1}} = 4k^{2}\frac{\cosh(x)}{\sinh(x)^{2k + 1}} + (2k + 2)(2k + 1)\frac{\cosh(x)}{\sinh(x)^{2k + 3}}
\end{align*}
we can derive the recurrence relation for the $c_{n, k}$'s
\begin{equation}
c_{n, 1} = \frac{1}{2}4^{n}, \quad c_{n, k} = 2k(2k - 1)c_{n - 1, k - 1} + 4k^{2}c_{n - 1, k} \;\; \mathrm{for} \;\; 2 \geq k \geq n - 1 \quad \mathrm{and} \;\; c_{n, n} = (2n)!.
\label{c_n_k_recurrence}
\end{equation}
Now the choice 
\begin{equation}
c_{n, k} = \frac{2}{4^{k}}\sum_{j = 1}^{k}(-1)^{k - j}\binom{2k}{k - j}(2j)^{2n} \label{definition_of_c_n_k}
\end{equation}
fulfills the first two of the latter equations even for $n = 0$, which can be readily checked by a small calculation using
\begin{align*}
\binom{2k - 2}{k - 1 - j} = \frac{k^2 - j^2}{2k(2k - 1)}\binom{2k}{k - j}.
\end{align*}

The trivial identity $2\sum_{j = 0}^{k - 1}(-1)^{j}\binom{2k}{j} = (-1)^{k + 1}\binom{2k}{k}$ directly yields $c_{0, k} = (-1)^{k + 1}\frac{1}{4^{k}}\binom{2k}{k}$ for $k \geq 1$. This in turn implies $c_{1, k} = 0$ for $k \geq 2$, which is due to $ c_{1, k} = 2k(2k - 1)c_{0, k - 1} + 4k^{2}c_{0, k}$. Then the same argument gives $c_{n, k} = 0$ for $k > n$. This gives $c_{n, n} = (2n)!$. Thus the claim follows by induction over $n$.
\hfill $\square$
\end{prf}

\lmm{The $c_{n, k}$'s fulfilling \eqref{c_n_k_recurrence} also fulfill for $n \geq k$
\begin{equation}
\sum_{i = k}^{n}\binom{k}{i - k}c_{n, i} = 2^{2(n - k)}c_{n, k}. \label{c_n_k_binomial_recurrence}
\end{equation}
}

\begin{prf}
We set $a_{n, k} :=\sum_{i = k}^{n}\binom{k}{i - k}c_{n, i}$ and $b_{n, k} := 2^{2(n - k)}c_{n, k}$. The strategy of the proof is to show that both $a_{n, k}$ and $b_{n, k}$ suffice
\begin{equation}
a_{n, k} = 2k(2k - 1)a_{n - 1, k - 1} + 16k^2 a_{n - 1, k}. \label{a_n_k_recurrence} 
\end{equation} 
Then since we have $a_{n, 1} = 2^{4n - 3} = b_{n, 1}$ and $a_{n, k} = b_{n, k} = 0$ for $k > n$ this gives $a_{n, k} = b_{n, k}$ for all $k \geq 1$ as claimed. Now \eqref{a_n_k_recurrence} is easily confirmed for the $b_{n, k}$'s. For the $a_{n, k}$'s we note that for $i \geq k$ 
\begin{align*}
(2i + 2)(2i + 1)\binom{k}{i - k + 1} - 2k(2k - 1)\binom{k - 1}{i - k + 1} + 4(i^2 - 4k^2)\binom{k}{i - k} = 0
\end{align*}
holds. Therefore we have
\begin{align*}
& \sum_{i = k}^{n - 1}(2i + 2)(2i + 1)\binom{k}{i - k + 1}c_{n - 1, i} + \sum_{i = k}^{n - 1}4 i^2\binom{k}{i - k}c_{n - 1, i} \\
= \;\; & 2k(2k - 1)\sum_{i = k}^{n - 1}\binom{k - 1}{i - k + 1}c_{n - 1, i} + 16k^2\sum_{i = k}^{n - 1}\binom{k}{i - k}c_{n - 1, i}.
\end{align*}
An index shift in the first sum of the left hand side and adding $2k(2k - 1)c_{n - 1, k - 1}$ on both sides yield
\begin{align*}
a_{n, k} = & \sum_{i = k}^{n}\binom{k}{i - k}\left(2i(2i - 1)c_{n - 1, i - 1} + 4 i^2c_{n - 1, i}\right) \\
= \;\; & 2k(2k - 1)\sum_{i = k - 1}^{n - 1}\binom{k - 1}{i - k + 1}c_{n - 1, i} + 16k^2\sum_{i = k}^{n - 1}\binom{k}{i - k}c_{n - 1, i} \\
= \;\; & 2k(2k - 1)a_{n - 1, k - 1} + 16k^2a_{n - 1, k}.
\end{align*}
\hfill $\square$
\end{prf}

\lmm{Let the matrix $U = (u_{n, k})$ be defined by 
\begin{equation}
u_{n, k} :=  \frac{(-1)^n}{(2n)!}c_{n, k}, \label{definition_of_U}
\end{equation}
with the $c_{n, k}$'s from \eqref{definition_of_c_n_k}, i.e. the $u_{n, k}$'s satisfy
\begin{equation}
u_{n, k} = - \frac{2k(2k - 1)}{2n(2n - 1)}u_{n - 1, k - 1} - \frac{4k^2}{2n(2n - 1)}u_{n - 1, k}. \label{u_recurrence}
\end{equation}
Let $V = (v_{n, k})$ denote the inverse of $U$. Then the $v_{n, k}$'s satisfy the recurrence relation
\begin{equation}
v_{n, k} = - \frac{2k(2k - 1)}{2n(2n - 1)}v_{n - 1, k - 1} - \frac{4(n - 1)^2}{2n(2n - 1)}v_{n - 1, k}. \label{v_recurrence}
\end{equation} 
}

\begin{prf}
Note that $U$ and $V$ are lower triangular matrices. The proof follows by induction over the row index $n$ for $V$. We fix an $s \in \{1, ..., n\}$ and compute
\begin{align*}
\sum_{k = s}^{n}v_{n, k}u_{k, s} = & - \sum_{k = s}^{n}\frac{2k(2k - 1)}{2n(2n - 1)}v_{n - 1, k - 1}u_{k, s} - \sum_{k = s}^{n}\frac{4(n - 1)^2}{2n(2n - 1)}v_{n - 1, k}u_{k, s} \\
= \;\; & - \frac{4(n - 1)^2}{2n(2n - 1)}\sum_{k = s}^{n - 1}v_{n - 1, k}u_{k, s} + \frac{2s(2s - 1)}{2n(2n - 1)}\sum_{k = s}^{n}v_{n - 1, k - 1}u_{k - 1, s - 1} \\ 
& + \frac{4s^2}{2n(2n - 1)}\sum_{k = s + 1}^{n}v_{n - 1, k - 1}u_{k - 1, s}
\end{align*}
All of the last three sums vanish for $s < n - 1$ by the induction hypothesis. For $s = n - 1$ we get 
\begin{align*}
\sum_{k = n - 1}^{n}v_{n, k}u_{k, n - 1} = - \frac{4(n - 1)^2}{2n(2n - 1)} + 0 +\frac{4(n - 1)^2}{2n(2n - 1)} = 0
\end{align*}
and for $s = n$ we obtain
\begin{align*}
\sum_{k = n}^{n}v_{n, k}u_{k, n} = 0 + 1 + 0.
\end{align*}
\hfill $\square$
\end{prf}

\prp{Let the numbers $h_{k, n}$ be given by
\begin{equation}
h_{1, n} = 1 \quad \mathrm{and} \quad h_{k, n} = \sum_{j = k - 1}^{n - 1}\frac{1}{j^2}h_{k - 1, j} \;\; \mathrm{for} \;\; k \geq 2. \label{definition_of_h_k_n}  
\end{equation}
Then the entries $v_{n, k}$ satisfying \eqref{v_recurrence} are given by
\begin{equation}
v_{n, k} = (-1)^{n}(2k)!\frac{2^{2(n - k)}}{n^2\binom{2n}{n}}h_{k, n}. \label{definition_of_V}
\end{equation}
}

\begin{prf}
Plugging the $v_{n, k}$'s as given in \eqref{definition_of_V} into \eqref{v_recurrence} gives a relation equivalent to 
\begin{align*}
h_{k, n} = \frac{1}{(n - 1)^2}h_{k - 1, n - 1} + h_{k, n - 1}.
\end{align*}
This in turn is equivalent to the definition \eqref{definition_of_h_k_n}. Furthermore we have
\begin{align*}
v_{n, n} =  (-1)^{n}(2n)!\frac{1}{n^2\binom{2n}{n}}h_{n, n} = (-1)^{n}((n - 1)!)^{2}\frac{1}{((n - 1)!)^{2}} = (-1)^{n}.
\end{align*}
\hfill $\square$
\end{prf}

\begin{prp}
We define the lower triangular matrix $L = (l_{n, k}) \in \mathbb{R}^{N \times N}$ by
\begin{align*}
l_{n, k} = \;\; & 2^{2k + 1}\binom{k}{n - k}, & & \mathrm{for} \;\; k \leq n \leq \min(2k, N) \\ 
l_{n, k} = \;\; & 0, & & \mathrm{otherwise}.
\end{align*}
and the diagonal matrix $D = (d_{n, k})$ by $d_{i, i} = 2^{2i + 1}$. Then we have for $U$ and $V$ defined in \eqref{definition_of_U} and \eqref{definition_of_V} respectively the matrix relations
\begin{align*}
UL = DU \quad \mathrm{and} \quad LV = VD.
\end{align*}
\end{prp}

\begin{prf}
Remembering $u_{n, k} = \frac{(-1)^n}{(2n)!}c_{n, k}$ we have
\begin{align*}
(UL)_{n, k} = \;\; & 2^{2k + 1}\frac{(-1)^n}{(2n)!}\sum_{i = k}^{n}\binom{k}{i - k}c_{n, i} = 2^{2k + 1}\frac{(-1)^n}{(2n)!}2^{2(n - k)}c_{n, k} = 2^{2n + 1}u_{n, k}.
\end{align*}
Here we used \eqref{c_n_k_binomial_recurrence} in the penultimate step.
\hfill $\square$
\end{prf}

\crl{The matrix equation $LV = VD$ can be written out as
\begin{align}
\sum_{k = 0}^{M - 1}2^{2(M + k) + 1}\frac{M - k}{2k + 1}\binom{M + k}{2k}v_{M + k, j} & = 2^{2j + 1}v_{2M - 1, j}  \label{V_odd_row_of_L} \\
\mathrm{and} \quad \sum_{k = 0}^{M}2^{2(M + k) + 1}\binom{M + k}{2k}v_{M + k, j} & = 2^{2j + 1}v_{2M, j} \label{V_even_row_of_L}
\end{align}
Note that the coefficients occurring above are the same as in \eqref{eq_coth_odd} and \eqref{eq_coth_even}.
}

\section{Main Results}

\subsection{A Limit Identity for $\lim_{\alpha \to 0_{+}}\sum_{n = 1}^{\infty}\frac{1}{n}\frac{\sinh(\alpha n)}{\cosh(\alpha n)^{2N + 1}}$}

Our starting point is Ramanujan's famous formula for $\zeta(2n + 1)$, cf. \cite{B74}.

\trm{Let $B_{r}$, $r \geq 0$ denote the $r$-th Bernoulli number. If $\alpha$ and $\beta$ are positive numbers sucht that $\alpha\beta = \pi^2$, and if $n$ is a positive integer, then
\begin{align*}
& \frac{1}{(4\alpha)^{n}}\left(\frac{1}{2}\zeta(2n + 1) + \sum_{m = 1}^{\infty}\frac{1}{m^{2n + 1}(\exp(2m\alpha) - 1)}\right) \\ 
- \;\; & \frac{1}{(- 4\beta)^{n}}\left(\frac{1}{2}\zeta(2n + 1) + \sum_{m = 1}^{\infty}\frac{1}{m^{2n + 1}(\exp(2m\beta) - 1)}\right) \\
= \;\; & \sum_{k = 0}^{n + 1}(-1)^{k - 1}\frac{B_{2k}}{(2k)!}\frac{B_{2n - 2k + 2}}{(2n - 2k + 2)!}\alpha^{n - k + 1}\beta^{k}.
\end{align*}
}

Using Euler's classical result for $n \in \mathbb{N}$
\begin{align*}
\zeta(2n) = \frac{(-1)^{n - 1}B_{2n}}{2(2n)!}(2\pi)^{2n},
\end{align*}
we rewrite Ramanujan's formula in a more convenient form, which is 

\crl{Let $M \in \mathbb{N}_{0}$. Then we have for $\alpha > 0$ and for $M \geq 1$
\begin{equation}
\begin{split}
& \pi\frac{1}{\alpha^{2M}}\sum_{n = 1}^{\infty}\frac{\coth(\alpha\pi n)}{n^{4M + 1}} - \pi\alpha^{2M}\sum_{n = 1}^{\infty}\frac{\coth(\frac{1}{\alpha}\pi n)}{n^{4M + 1}} \\
= \;\; & - \zeta(4M + 2)\left(\alpha^{2M + 1} - \frac{1}{\alpha^{2M + 1}}\right) + 2\sum_{j = 1}^{M}(-1)^{j}\zeta(2j)\zeta(4M + 2 - 2j)\left(\alpha^{2M + 1 - 2j} - \frac{1}{\alpha^{2M + 1 - 2j}}\right) \\
= \;\; & \zeta(4M + 2)\frac{1}{\alpha^{2M + 1}} - 2\sum_{k = 1}^{2M}(-1)^{k}\zeta(2k)\zeta(4M + 2 - 2k)\frac{1}{\alpha^{2(M - k) + 1}} - \zeta(4M + 2)\alpha^{2M + 1}, \label{eq_N_even}
\end{split}
\end{equation}
and for $M \geq 0$
\begin{equation}´\
\begin{split}
& \pi\frac{1}{\alpha^{2M + 1}}\sum_{n = 1}^{\infty}\frac{\coth(\alpha\pi n)}{n^{4M + 3}} + \pi\alpha^{2M + 1}\sum_{n = 1}^{\infty}\frac{\coth(\frac{1}{\alpha}\pi n)}{n^{4M + 3}} \\
= \;\; & \zeta(4M + 4)\left(\alpha^{2M + 2} + \frac{1}{\alpha^{2M + 2}}\right) - 2\sum_{j = 1}^{M}(-1)^{j}\zeta(2j)\zeta(4M + 4 - 2j)\left(\alpha^{2(M + 1 - j)} + \frac{1}{\alpha^{2(M + 1 - j)}}\right) \\
+ \;\; & 2(-1)^{M}\zeta(2M + 2)^{2} \\
= \;\; & \zeta(4M + 4)\frac{1}{\alpha^{2M + 2}} - 2\sum_{k = 1}^{2M + 1}(-1)^{k}\zeta(2k)\zeta(4M + 4 - 2k)\frac{1}{\alpha^{2(M - k + 1)}} + \zeta(4M + 4)\alpha^{2M + 2}. \label{eq_N_odd}
\end{split}
\end{equation}
}

\trm{Let $N \in \mathbb{N}$ and $h_{k, n}$ as defined in \eqref{definition_of_h_k_n}. Then we have
\begin{equation}
\lim_{\alpha \to 0_{+}}\sum_{n = 1}^{\infty}\frac{1}{n}\frac{\sinh(\alpha n)}{\cosh(\alpha n)^{2N + 1}} = \frac{1}{N^2\binom{2N}{N}}\sum_{k = 1}^{N}(2k)!\left(2^{2N + 1} - 2^{2(N - k)}\right)h_{k, N}\frac{\zeta(2k + 1)}{\pi^{2k}}. \label{limit_equation}
\end{equation}
}

\begin{prf}
We divide both sides of \eqref{eq_N_odd} by $\alpha^{2M + 1}$ and both sides of \eqref{eq_N_even} by $\alpha^{2M}$. Then we apply the operator $\alpha^2\frac{d}{d\alpha}$ on \eqref{eq_N_odd} $(4M + 2)$ times and on \eqref{eq_N_even} $4M$ times. This results in the linear system
\begin{align*}
\pi^{4M + 3}\sum_{k = 1}^{2M + 1}c_{2M + 1, k}\sum_{n = 1}^{\infty}\frac{1}{n}\frac{\cosh(\frac{1}{\alpha}\pi n)}{\sinh(\frac{1}{\alpha}\pi n)^{2k + 1}} = \;\;  & - \pi(4M + 2)!\sum_{n = 1}^{\infty}\frac{\coth(\alpha\pi n)}{n^{4M + 3}} \\
 & + \; \{\mathrm{sums \; involving \; \coth \; vanishing \; for \; \alpha \rightarrow \infty}\} \\
 & + \; (4M + 3)!\zeta(4M + 4)\frac{1}{\alpha} + (4M + 2)!\zeta(4M + 4)\alpha^{4M + 3}, 
\end{align*}
for $M \geq 0$ and
\begin{align*}
- \pi^{4M + 1}\sum_{k = 1}^{2M}c_{2M, k}\sum_{n = 1}^{\infty}\frac{1}{n}\frac{\cosh(\frac{1}{\alpha}\pi n)}{\sinh(\frac{1}{\alpha}\pi n)^{2k + 1}} = \;\; & - \pi(4M)!\sum_{n = 1}^{\infty}\frac{\coth(\alpha\pi n)}{n^{4M + 1}} \\
 & + \; \{\mathrm{sums \; involving \; \coth \; vanishing \; for \; \alpha \rightarrow \infty}\} \\
 & + \; (4M + 1)!\zeta(4M + 2)\frac{1}{\alpha} - (4M)!\zeta(4M + 2)\alpha^{4M + 1}, \;\:\:
\end{align*}
for $M \geq 1$. Here we use the $c_{n, k}$'s from \eqref{definition_of_c_n_k}. Select $N \in \mathbb{N}$. Then solving for $\sum_{n = 1}^{\infty}\frac{1}{n}\frac{\cosh(\frac{1}{\alpha}\pi n)}{\sinh(\frac{1}{\alpha}\pi n)^{2N + 1}}$ gives
\begin{align*}
\sum_{n = 1}^{\infty}\frac{1}{n}\frac{\cosh(\frac{1}{\alpha}\pi n)}{\sinh(\frac{1}{\alpha}\pi n)^{2N + 1}} = \;\; & \sum_{k = 1}^{N}v_{N, k}\frac{1}{\pi^{2k}}\sum_{n = 1}^{\infty}\frac{\coth(\alpha\pi n)}{n^{2k + 1}} \\
& + \; \{\mathrm{sums \; involving \; \coth \; vanishing \; for \; \alpha \rightarrow \infty}\} \\
& - \sum_{k = 1}^{N}v_{N, k}(2k + 1)\frac{\zeta(2k + 2)}{\pi^{2k + 1}}\frac{1}{\alpha} + \sum_{k = 1}^{N}v_{N, k}(-1)^{k}\frac{\zeta(2k + 2)}{\pi^{2k + 1}}\alpha^{2k + 1},
\end{align*}
with the $v_{n, k}$'s defined in \eqref{definition_of_V}. Let us first assume that $N$ is odd, i.e. $N = 2M - 1$. Then plugging \eqref{eq_coth_odd} into the left hand sind of above equation gives
\begin{align*}
\sum_{n = 1}^{\infty}\frac{1}{n}\frac{\sinh(\frac{1}{\alpha}\pi n)}{\cosh(\frac{1}{\alpha}\pi n)^{4M - 1}}  = \;\; & - \sum_{k = 0}^{M - 1}2^{2(M + k) + 1}\frac{M - k}{2k + 1}\binom{M + k}{2k}\sum_{n = 1}^{\infty}\frac{1}{n}\frac{\cosh(2\frac{1}{\alpha}\pi n)}{\sinh(2\frac{1}{\alpha}\pi n)^{2(M + k) + 1}} \\
& + \sum_{k = 1}^{2M - 1}v_{2M - 1, k}\frac{1}{\pi^{2k}}\sum_{n = 1}^{\infty}\frac{\coth(\alpha\pi n)}{n^{2k + 1}} \\
& + \; \{\mathrm{sums \; involving \; \coth \; vanishing \; for \; \alpha \rightarrow \infty}\} \\
& - \sum_{k = 1}^{2M - 1}v_{2M - 1, k}(2k + 1)\frac{\zeta(2k + 2)}{\pi^{2k + 1}}\frac{1}{\alpha} + \sum_{k = 1}^{2M - 1}v_{2M - 1, k}(-1)^{k}\frac{\zeta(2k + 2)}{\pi^{2k + 1}}\alpha^{2k + 1} \\
= \;\; & - \sum_{k = 0}^{M - 1}2^{2(M + k) + 1}\frac{M - k}{2k + 1}\binom{M + k}{2k}\sum_{s = 1}^{M + k}v_{M + k,s}\frac{1}{\pi^{2s}}\sum_{n = 1}^{\infty}\frac{\coth(\frac{1}{2}\alpha\pi n)}{n^{2s + 1}}\\
& + \sum_{k = 1}^{2M - 1}v_{2M - 1, k}\frac{1}{\pi^{2k}}\sum_{n = 1}^{\infty}\frac{\coth(\alpha\pi n)}{n^{2k + 1}} \\
& + \; \{\mathrm{sums \; involving \; \coth \; vanishing \; for \; \alpha \rightarrow \infty}\} \\
& - \sum_{k = 0}^{M - 1}2^{2(M + k) + 1}\frac{M - k}{2k + 1}\binom{M + k}{2k}\sum_{s = 1}^{M + k}\bigg( - v_{M + k, s}(2s + 1)\frac{\zeta(2s + 2)}{\pi^{2s + 1}}\frac{2}{\alpha} \\
& + v_{M + k, s}(-1)^{s}\frac{\zeta(2s + 2)}{\pi^{2s + 1}}\left(\textstyle{\frac{1}{2}}\alpha\right)^{2s + 1}\bigg) \\
& - \sum_{k = 1}^{2M - 1}v_{2M - 1, k}(2k + 1)\frac{\zeta(2k + 2)}{\pi^{2k + 1}}\frac{1}{\alpha} + \sum_{k = 1}^{2M - 1}v_{2M - 1, k}(-1)^{k}\frac{\zeta(2k + 2)}{\pi^{2k + 1}}\alpha^{2k + 1} \\
= \;\; & - \sum_{s = 1}^{M}\sum_{k = 0}^{M - 1}2^{2(M + k) + 1}\frac{M - k}{2k + 1}\binom{M + k}{2k}v_{M + k,s}\frac{1}{\pi^{2s}}\sum_{n = 1}^{\infty}\frac{\coth(\frac{1}{2}\alpha\pi n)}{n^{2s + 1}} \\
& - \sum_{s = 1}^{M - 1}\sum_{k = s}^{M - 1}2^{2(M + k) + 1}\frac{M - k}{2k + 1}\binom{M + k}{2k}v_{M + k, M + s}\frac{1}{\pi^{2(M + s)}}\sum_{n = 1}^{\infty}\frac{\coth(\frac{1}{2}\alpha\pi n)}{n^{2(M + s) + 1}} \\
& + \sum_{k = 1}^{2M - 1}v_{2M - 1, k}\frac{1}{\pi^{2k}}\sum_{n = 1}^{\infty}\frac{\coth(\alpha\pi n)}{n^{2k + 1}} \\
& + \; \{\mathrm{sums \; involving \; \coth \; vanishing \; for \; \alpha \rightarrow \infty}\} \\
& + c_{2M - 1}\frac{1}{\alpha} + P_{2M - 1}(\alpha), 
\end{align*}
where $c_{2M - 1}$ denotes the coefficient of $\frac{1}{\alpha}$ and $P_{2M - 1}(\alpha)$ the polynomial of degree $4M - 1$ occurring right after $ c_{2M - 1}\frac{1}{\alpha}$. Now \eqref{V_odd_row_of_L} and $v_{n, k} = 0$ for $k > n$ yield
\begin{equation}
\begin{split}
\sum_{n = 1}^{\infty}\frac{1}{n}\frac{\sinh(\frac{1}{\alpha}\pi n)}{\cosh(\frac{1}{\alpha}\pi n)^{4M - 1}}  = \;\; & \sum_{k = 1}^{2M - 1}v_{2M - 1, k}\frac{1}{\pi^{2k}}\left(\sum_{n = 1}^{\infty}\frac{\coth(\alpha\pi n)}{n^{2k + 1}} - 2^{2k + 1}\sum_{n = 1}^{\infty}\frac{\coth(\frac{1}{2}\alpha\pi n)}{n^{2k + 1}}\right) \\
& + \; \{\mathrm{sums \; involving \; \coth \; vanishing \; for \; \alpha \rightarrow \infty}\} \\
& + c_{2M - 1}\frac{1}{\alpha} + P_{2M - 1}(\alpha). \label{odd_intermediate_result}
\end{split}
\end{equation}
For $M = 1$ we can with a little effort directly compute 
\begin{align*}
\lim_{\alpha \to \infty}\sum_{n = 1}^{\infty}\frac{1}{n}\frac{\sinh(\frac{1}{\alpha}\pi n)}{\cosh(\frac{1}{\alpha}\pi n)^{3}} = \frac{7}{\pi^2}\zeta(3).
\end{align*}
Since $\frac{\sinh(\frac{1}{\alpha}\pi n)}{\cosh(\frac{1}{\alpha}\pi n)^{2M + 1}} \leq \frac{\sinh(\frac{1}{\alpha}\pi n)}{\cosh(\frac{1}{\alpha}\pi n)^{3}}$ for all $M \geq 1$ this gives
\begin{align*}
0 < \sum_{n = 1}^{\infty}\frac{1}{n}\frac{\sinh(\frac{1}{\alpha}\pi n)}{\cosh(\frac{1}{\alpha}\pi n)^{2M + 1}} < \frac{7}{\pi^2}\zeta(3) + \epsilon
\end{align*}  
for a given $\epsilon$ and $\alpha$ sufficiently large. All other terms than $P_{2M - 1}(\alpha)$ in the right hand side of \eqref{odd_intermediate_result} are bounded. Since $P_{2M - 1}(\alpha)$ does not have a constant term it follows with above result $P_{2M - 1}(\alpha) \equiv 0$. Inserting \eqref{definition_of_V} into \eqref{odd_intermediate_result} and the fact that $\lim_{\alpha \to \infty}\coth(\alpha\pi n) = 1$ holds for all $n \in \mathbb{N}$ proves the claim for $N = 2M - 1$. The proof for $N = 2M$ follows completely analogously using \eqref{eq_coth_even} and \eqref{V_even_row_of_L}.

\hfill $\square$
\end{prf}

\subsection{Applications}

\crl{We have the identity
\begin{equation}
14\sum_{n = 1}^{\infty}\frac{n2^{n}}{(n + 1)\binom{2n + 2}{n + 1}}\frac{\zeta(3)}{\pi^{2}} + \sum_{j = 2}^{\infty}\left(\sum_{n = j - 1}^{\infty}\frac{n2^{n}}{(n + 1)\binom{2n + 2}{n + 1}}h_{j, n + 1}\right)4(2j)!\left(2 - \frac{1}{2^{2j}}\right)\frac{\zeta(2j + 1)}{\pi^{2j}} = 28\frac{\zeta(3)}{\pi^2}. \label{odd_zeta_recurrence1}
\end{equation}
}

\begin{prf}
Differentiating
\begin{equation}
\tanh(x)\sum_{k = 0}^{\infty}\left(\frac{1}{2\cosh(x)^{2}}\right)^{k} = \tanh(2x) \label{id_tanh}
\end{equation}
twice gives
\begin{equation}
\sum_{k = 1}^{\infty}\frac{k(k + 1)}{2^{k}}\frac{\sinh(x)}{\cosh(x)^{2k + 3}} = 4\frac{\sinh(2x)}{\cosh(2x)^{3}}. \label{id_from_tanh}
\end{equation}
Inserting \eqref{limit_equation}, taking the limit $\alpha \to 0_{+}$ and rearranging give the result.
\hfill $\square$
\end{prf}

\rem{More identities like \eqref{odd_zeta_recurrence1} can be obtained from differentiating \eqref{id_from_tanh} $2N$ times.}

\crl{Choose $K \in \mathbb{N}_{0}$ and $N \in \mathbb{N}$. Then the limits
\begin{equation}
\lim_{\alpha \to 0_{+}}\frac{1}{\alpha^{K}}\sum_{n = 1}^{\infty}\frac{1}{n^{K + 1}}\frac{\sinh(\alpha n)^{1 + K}}{\cosh(\alpha n)^{2M + 1 + K}} \label{K_equation}
\end{equation}
exist and are a finite linear combination of $\frac{\zeta(3)}{\pi^2}$, ..., $\frac{\zeta(2N + 2K + 1)}{\pi^{2N + 2K}}$ with rational coefficients.
}

\begin{prf}
For $K = 0$ we the claim becomes \eqref{limit_equation}. For $K > 0$ l\textsc{\char13}H\^{o}pital's rule yields
\begin{align*}
\lim_{\alpha \to 0_{+}}\frac{1}{\alpha^{K}}\sum_{n = 1}^{\infty}\frac{1}{n^{K + 1}}\frac{\sinh(\alpha n)^{1 + K}}{\cosh(\alpha n)^{2M + 1 + K}} = \;\; & -\frac{2N}{K}\lim_{\alpha \to 0_{+}}\frac{1}{\alpha^{K - 1}}\sum_{n = 1}^{\infty}\frac{1}{n^{K}}\frac{\sinh(\alpha n)^{1 + K - 1}}{\cosh(\alpha n)^{2M + 1 + K - 1}} \\
& + \frac{K + 2N + 1}{K}\lim_{\alpha \to 0_{+}}\frac{1}{\alpha^{K - 1}}\sum_{n = 1}^{\infty}\frac{1}{n^{K}}\frac{\sinh(\alpha n)^{1 + K - 1}}{\cosh(\alpha n)^{2M + 3 + K - 1}}.
\end{align*}
Iterating this finding $K$ times and \eqref{limit_equation} give the result.
\hfill $\square$
\end{prf}

\rem{Using \eqref{K_equation} one can obtain for $K \in \mathbb{N}$ more identities like \eqref{odd_zeta_recurrence1} from \eqref{id_tanh} and \eqref{id_from_tanh} by sending $\alpha \to 0_{+}$ in
\begin{align*}
& \frac{1}{\alpha^{K}}\sum_{n = 1}^{\infty}\frac{1}{n^{1 + K}}\left(\tanh(\alpha n)\sum_{k = 0}^{\infty}\left(\frac{1}{2\cosh(\alpha n)^{2}}\right)^{k}\right)^{K}\sum_{k = 1}^{\infty}\frac{k(k + 1)}{2^{k}}\frac{\sinh(\alpha n)}{\cosh(\alpha n)^{2k + 3}}\\ 
= \;\; & 4\frac{1}{\alpha^{K}}\sum_{n = 1}^{\infty}\frac{1}{n^{1 + K}}\tanh(2\alpha n)^{K}\frac{\sinh(2\alpha n)}{\cosh(2\alpha n)^{3}}.
\end{align*}
}

\bibliographystyle{ieeetr}
\bibliography{literature}

\end{document}